\documentclass{article}
\usepackage{amssymb} 
\title{{\bf Random homogenization  of $p$-Laplacian with obstacles  in perforated domain }
\vspace*{0.2cm}
\author{\small  Lan Tang\\
\small Department of  Mathematics\\
\small The University of Texas at Austin\\
\small Austin, TX 78712 \\
\small }
\date{}}

\begin{document}
\maketitle
\begin{minipage}{14 true cm}{\bf Abstract:} \  {In this paper,
we will study the homogenization  of $p$-Laplacian with obstacles in perforated domain, where the holes are
periodically distributed and have random size. And we also assume that the $p$-capacity of each hole is stationary ergodic.}\\

\end{minipage}

\section{ Introduction} 
 
Let $D\subset{\bf R}^n$ be a bounded domain and $(\Omega, \mathcal {F}, \mathcal {P})$ be
a given probability space. For each $\omega \in \Omega$ and
$\epsilon >0$, we denote by $T_{\epsilon}(\omega)$ the set of holes on $D$. 
Our main purpose  is to study the asymptotic behavior as $\epsilon\rightarrow0$ 
of the solution $u^\epsilon $ of the following variational problem:
$$\displaystyle{\min}\big\{\int_{D}\frac{1}{p}{|\nabla u|^p}dx-\int_{D}fudx: u \in W^{1,p}_{0}(D),\  u\geq 0 \ \mbox{a.\ e.} \ \mbox{in} \ T_{\epsilon}(\omega) \big\}
$$
where $f$ is some measurable and bounded function.

This is a classical homogenization problem  and the asymptotic
 behavior of  $u^\epsilon $ strongly depends on the properties of
 $T_{\epsilon}(\omega)$. This type of problems were first studied by
 L. Carbone and F. Colombini [CC] in periodic settings and then in more general frameworks by
 E. De Giorgi, G. Dal Maso and P. Longo[DDL],G. Dal Maso and P. Longo [DL] and G. Dal Maso[D]. And D. Cioranescu and F. Murat [CM1,2] studied the special case $p=2$ (Laplacian) for
 periodic settings. For more general structure, [AB] used $\Gamma$-convergence to study the periodic homogenization, which generalized the result of D. Cioranescu and F. Murat [CM1,2].
 For the stationary ergodic settings, Caffarelli and Mellet [CM] studied the
 case $p=2$ and Focardi [F] used $\Gamma$-convergence method to study fractional obstacle problems. 

In this article , we consider the  case for the $1<p\leq n$. 
We assume that  $T_{\epsilon}(\omega)$ satisfies the following:
$$T_{\epsilon}(\omega)=\bigcup_{k\in{\bf Z}^n}B_{a^{\epsilon}(k,\omega)}(\epsilon k) $$
and the p-capacity (see [MZ]) of  each ball $B_{a^{\epsilon}(k,\omega)}(\epsilon k)$ satisfies :
$$\mbox{cap}_{p}(B_{a^{\epsilon}(k,\omega)}(\epsilon k))=\gamma(k, \omega)\epsilon^{n}$$
where $\gamma: {\bf Z}^n \times \Omega \mapsto [0,+\infty) $ is 
stationary ergodic: there exists a family of measure-preserving transformations $\tau_{k}: \Omega \mapsto \Omega$ satisfying
$$\gamma(k+k^{'}, \omega)=\gamma(k, \tau_{k^{'}}\omega),  \ \forall \ k,\  k^{'} \in {\bf Z}^n \ \mbox{and} \  \omega \in \Omega,$$
and such that if $A \subset\Omega$ and $\tau_{k}A=A$ for all $k\in
{\bf Z}^n$, then $P(A)=1$ or $P(A)=0$. And we also assume that $\gamma: {\bf Z}^n \times \Omega \mapsto [0,+\infty) $ is bounded.

Thus
\[ a^{\epsilon}(k, \omega)=\left 
\{\begin{array}{ll}
(\frac{\gamma(k, \omega)}{n\omega_{n}})^{\frac{1}{n-p}}(\frac{n-p}{p-1})^{\frac{1-p}{n-p}}{\epsilon}^{\frac{n}{n-p}} &  \  \quad  \mbox{if} \qquad 1<p<n\\
\exp(-(\frac{\gamma(k, \omega)}{n\omega_{n}})^{\frac{-1}{n-1}}{\epsilon}^{-\frac{n}{n-1}})&  \  \quad \mbox{if} \qquad p=n
\end{array}
\right. \]
Obviously, $a^{\epsilon}: {\bf Z}^n \times \Omega \mapsto [0,+\infty) $ is also
stationary ergodic and bounded.
 
In the following, we firstly consider the variational problem :
$$\inf_{v\in K_{\epsilon}}{\mathcal{F}}(v)$$
where  $ \displaystyle{\mathcal{F}}(v)=\int_{D}{\frac{1}{p}|\nabla{v}|^p-fv}dx$
and  $\displaystyle K_{\epsilon}=\{v\in W_{0}^{1, p}(D): v\geq0 \ \mbox{a.\ e.}\ \mbox{on} \
 T_{\epsilon}\}$.

Let $u^{\epsilon}$ be the solution of such a variational problem, i.e.
$${\mathcal{F}}(u^{\epsilon})=\inf_{v\in K_{\epsilon}}{\mathcal{F}}(v)$$.

Obviously, $\{u^{\epsilon}\}$ is bounded in
 $W_{0}^{1, p}(D)$, then we can choose a subsequence of $\{u^{\epsilon}\}$ (we still denote by $u^{\epsilon}$) such that
 $$u^{\epsilon}\rightharpoonup {u}^{0} \ \ \mbox{in} \ \ W_{0}^{1, p}(D).$$
Our main purpose is to determine the  variational functional ${\mathcal{F}}_{0}$ 
such that for almost surely $\omega\in\Omega$,\ 
$${\mathcal{F}}_{0}({u}^{0})=\inf_{v\in
 W_{0}^{1, p}(D)}{{\mathcal F}_{0}}(v).$$

Now we state our main results: \\

{\bf  Theorem 1.1.}\textit {\quad Let $1< p\leq n$ and the stochastic process $\gamma(k, \omega): Z^n \times \Omega \rightarrow [0, \infty)$ is bounded above by some universal positive constant 
then there exits a nonnegative real number $\alpha_{0}$
such that when $\epsilon$ goes to zero, the solution $u^{\epsilon}(x,\omega)$
of
$$\displaystyle \min\{\int_{{\bf R}^n}\frac{1}{p}|\nabla{v}|^{p}-fvdx: \  u\in W_{0}^{1,p}(D),\  u\geq 0 \  \mbox{a.\ e.}\  \mbox{in}\  T_{\epsilon}(\omega)\}
$$\\
converges weakly in $W^{1,p}(D)$ and almost surely $\omega \in
\Omega$ to the solution $u_{0}$ of the following minimization problem:
$$\min\{ \int_{D}\frac{1}{p}|\nabla{v}|^{p}+\frac{1}{p}\alpha_{0}v_{-}^{p}-fvdx : \  \forall \  v\in W_{0}^{1,p}(D)\}$$
}

Next we consider the following variational inequality with oscillating obstacles:
$$\min\{\int_{D}\frac{1}{p}|\nabla v|^p -fv dx : \ v \in W^{1,p}_{0}(D) \ \mbox{and} \ \ v\geq\psi^{\epsilon}\}$$
where $\psi$ be a measurable function in $D$ and 
\[ \psi^{\epsilon}=\left 
\{\begin{array}{ll}
 \psi &  \  \quad  \mbox{in} \qquad D\setminus T_{\epsilon}\\
 0 &  \  \quad \mbox{on} \qquad T_{\epsilon}
\end{array}
\right. \]

Let us suppose that $h^{\epsilon}$ is the solution to the problem above, then $h^{\epsilon}$ is obviously bounded in $W^{1,p}_{0}(D)$.
Hence there is some function $h_{0}$ such that $h^{\epsilon}$ converges to $h_0$ weakly in $W^{1,p}_{0}(D)$. 
Then we have the following result:
\\
{ \bf Corollary 1.2.}\textit {\quad For $1<p\leq n$, if  when $\epsilon$ goes to zero, the solution $h^{\epsilon}(x,\omega)$
of
$$\displaystyle \min\{\int_{{\bf R}^n}\frac{1}{p}|\nabla{v}|^{p}-fvdx: \  v\in W_{0}^{1,p}(D),\  v\geq \psi^{\epsilon} \  \mbox{a.\ e.}\  \mbox{in}\ D
$$converges weakly in $W^{1,p}(D)$ and almost surely $\omega \in\Omega$ to the solution $h_{0}$ , then $h_0$ is the solution to the following variational problem:
$$\min\{\int_{D}\frac{1}{p}|\nabla{v}|^{p}+\frac{1}{p}\alpha_{0}v_{-}^{p}-fvdx: \ v\in W_{0}^{1,p}(D) \ \mbox{and} \ v\geq\psi \
 \mbox{a.}\ \mbox{e.}\ \mbox{in} \  D\} $$
where the constant $\alpha_{0}$ is the same constant as in Theorem 1.1.}

\section{ Proof of the Main Theorem} 
As in D. Cioranescu and F. Murat [CM1,2] and Caffarelli-Mellet [CM], 
the proof of Theorem 1.1 and Corollary 1.2 depend on the properties of some suitable correctors. 
More precisely, we need the following two lemmas:

{\bf  Key Lemma I } \textit{ \quad Assume that $T_{\epsilon}(\omega)$ satisfies the assumptions listed above. Then there
exist a nonnegative real number $\alpha_{0}$  and a function $w^{\epsilon}$  such that
\[ \left\{\begin{array}{lll}
\displaystyle{\triangle}_{p} w^{\epsilon}&= \alpha_{0} \; \quad \mbox{in} \qquad D_{\epsilon}(\omega) \nonumber\\
w^{\epsilon}(x,\omega)&=1  \qquad  \mbox{for}\; \quad x \in T_{\epsilon}(\omega) \nonumber\\
w^{\epsilon}(x,\omega)&= 0 \qquad \mbox{for}\; \quad x\in
{\partial D}\setminus T_{\epsilon}(\omega) \nonumber\\
w^{\epsilon}(\cdot,\omega)&\rightarrow 0 \ \ \ \mbox{weakly}\  \mbox{in}  \ W^{1,p}\  \nonumber
\end{array}
\right.\]
for a. s. $ \omega \in \Omega$ and $w^{\epsilon}$ also satisfies the following properties: \\
{\qquad \bf (a)}\ for  any $\phi \in \mathcal{D}(D)$ and $0<p'<p$,
$$\lim_{\epsilon\rightarrow0}\int_{D}|\nabla w^{\epsilon}|^{p'}\phi dx =0$$\\
\qquad {\bf (b)} for  any $\phi \in \mathcal{D}(D)$,
$$\displaystyle\lim_{\epsilon\rightarrow0}\int_{D}|\nabla{w^{\epsilon}}|^{p}\phi dx=\int_{D}\alpha_{0}\phi dx .$$\\
\qquad {\bf(c)} \ for any sequence $\{v^{\epsilon}\}\subset W^{1,p}_{0}(D)$ with the property: 
$v^\epsilon \rightharpoonup v $ in $W^{1,p}_{0}(D)$ and $v^\epsilon =0$ on $T_{\epsilon}$ and any $\phi \in \mathcal{D}(D)$, 
we have that $$\lim_{\epsilon \rightarrow 0}\int_{D}|\nabla w^{\epsilon}|^{p-2}{\nabla w^{\epsilon}}\cdot {\nabla v^{\epsilon}}\phi dx = -\alpha_{0}\int_{D}v\phi dx.$$}

And for  $u^\epsilon$ in Theorem 1.1, we have the following lower semi-continuous  property: 
( For the special case : $p=2$,  we refer the readers to Proposition 3.1 of [CM 1,2] )\\

{\bf Key Lemma II} \textit{\quad If the key lemma I holds and $u^\epsilon$ is the solution of 
$$\displaystyle \min\{\int_{{\bf R}^n}\frac{1}{p}|\nabla{v}|^{p}-fvdx: \  u\in W_{0}^{1,p}(D),\  u\geq 0 \  \mbox{a.\ e.}\  \mbox{in}\  T_{\epsilon}(\omega)\}
$$
then $$\liminf_{\epsilon\rightarrow0}{\mathcal{F}}({u^{\epsilon}})\geq {\mathcal{F}}_{0}(u^{0}),$$
where $u^{0}$ is the weak limit of $\{u^{\epsilon}\}$ in $W^{1,p}$ and $\displaystyle{\mathcal{F}}_{0}$ is defined as follows:  
$${\mathcal{F}}_{0}(v)=  \int_{D}\frac{1}{p}|\nabla{v}|^{p}+\frac{1}{p}\alpha_{0}v_{-}^{p}-fvdx,\ \forall \ v\in  W^{1,p}_{0}(D).$$}\\

Now let us show that how the two lemmas I and II imply our main results:

{ \bf Proof of  Theorem 1.1:}\quad Let $\phi\in C^{1}_{0}(D) $
such that $\phi_{-}\in C^{1}_{0}(D)$. Then
$${\mathcal{F}}(u^{\epsilon})\leq{\mathcal{F}}(\phi+{\phi}_{-}w^{\epsilon}).$$
\quad Next, we can estimate
${\mathcal{F}}(\phi+{\phi}_{-}w^{\epsilon})$ as follows:
\begin{eqnarray*}
{\mathcal{F}}(\phi+{\phi}_{-}w^{\epsilon})&=& \int_{D}\frac{1}{p}[ \ |\nabla{\phi}+\nabla{\phi_{-}}{w^{\epsilon}}+\nabla{w^{\epsilon}}\phi_{-}|^p \ ]dx\\
&-&\int_{D}[f\phi+f\phi_{-}w^{\epsilon}]dx.
\end{eqnarray*}

If $p$ is an integer, then
\begin{eqnarray*}
|\nabla{\phi}+\nabla{\phi_{-}}{w^{\epsilon}}+\nabla{w^{\epsilon}}\phi_{-}|^p  &\leq&\{|\nabla{\phi}+\nabla{\phi_{-}}{w^{\epsilon}}| \ + \ |\nabla{w^{\epsilon}}\phi_{-}|\}^{p} \\
&=&\sum_{k=0}^{p}C_{p}^k |\nabla{\phi}+\nabla{\phi_{-}}{w^{\epsilon}}|^{k}\cdot|\nabla{w^{\epsilon}}\phi_{-}|^{p-k}\\
&=&|\nabla{\phi}+\nabla{\phi_{-}}{w^{\epsilon}}|^{p} \ + \ |\nabla{w^{\epsilon}}\phi_{-}|^{p} \\
&+& \sum_{k=1}^{p-1}C_{p}^k |\nabla{\phi}+\nabla{\phi_{-}}{w^{\epsilon}}|^{k}\cdot|\nabla{w^{\epsilon}}\phi_{-}|^{p-k}
\end{eqnarray*}
By Key Lemma I, $w^\epsilon$ converges to zero weakly in $W^{1,p}(D)$, then
$\displaystyle w^\epsilon \rightarrow 0 \ \mbox{strongly} \ \mbox{in} \ L^{p}(D) \ \mbox{as} \ \epsilon \rightarrow 0.$
Thus 
$$\lim_{\epsilon\rightarrow0}\int_{D}|\nabla{\phi}+\nabla{\phi_{-}}{w^{\epsilon}}|^{p}dx =\int_{D}|\nabla{\phi}|^{p}dx$$

By Key Lemma I, we have
$$\lim_{\epsilon\rightarrow0}\int_{D}|\nabla{w^{\epsilon}}\phi_{-}|^{p}dx
dx=\int_{D}\alpha_{0}{\phi_{-}}^p  dx$$ 

For any $k$: $1\leq k\leq p-1$, we have that $$C_{p}^k |\nabla{\phi}+\nabla{\phi_{-}}{w^{\epsilon}}|^{k} \leq C  + C|{w^{\epsilon}}|^{k}$$
where $C$ is a  constant (not depending on $k$).

From Key Lemma I, we know that ($1\leq k\leq p-1$)
$$\lim_{\epsilon\rightarrow0}\int_{D}|\nabla{w^{\epsilon}}\phi_{-}|^{p-k}dx =0$$

And by H\"{o}lder inequality,

$$\int_{D}|w^\epsilon|^k |\nabla{w^{\epsilon}}\phi_{-}|^{p-k}dx \leq \{\int_{D}|w^\epsilon|^p dx\}^{\frac{k}{p}}\cdot \{\int_{D}|\nabla{w^{\epsilon}}\phi_{-}|^{p}dx\}^{\frac{p-k}{p}}$$

Hence $$\lim_{\epsilon\rightarrow0}\int_{D}|w^\epsilon|^k |\nabla{w^{\epsilon}}\phi_{-}|^{p-k}dx=0.$$

Therefore (for $p$ is an integer)

$$\limsup_{\epsilon\rightarrow0}\int_{D}|\nabla{\phi}+\nabla{\phi_{-}}{w^{\epsilon}}+\nabla{w^{\epsilon}}\phi_{-}|^p dx \leq \int_{D}|\nabla \phi|^p dx + \int_{D}\alpha_{0}{\phi_{-}}^{p}dx$$

If $p$ is not an integer, then we let $m$ be the integer part of $p$ ( Thus $0< p-m <1$). Hence
\begin{eqnarray*}
 |\nabla{\phi}+\nabla{\phi_{-}}{w^{\epsilon}}+\nabla{w^{\epsilon}}\phi_{-}|^p  &\leq& \{ |\nabla{\phi}+\nabla{\phi_{-}}{w^{\epsilon}}|\ + \ |\nabla{w^{\epsilon}}\phi_{-}|\}^{p-m}\\
 &\times& \{|\nabla{\phi}+\nabla{\phi_{-}}{w^{\epsilon}}| \ + \ |\nabla{w^{\epsilon}}\phi_{-}|\}^{m} \\
&\leq& \{ |\nabla{\phi}+\nabla{\phi_{-}}{w^{\epsilon}}|^{p-m}\ + \ |\nabla{w^{\epsilon}}\phi_{-}|^{p-m}\}\\
&\times& \{|\nabla{\phi}+\nabla{\phi_{-}}{w^{\epsilon}}|\ + \ |\nabla{w^{\epsilon}}\phi_{-}|\}^{m}\\
&=& |\nabla{\phi}+\nabla{\phi_{-}}{w^{\epsilon}}|^p \ +  \ |\nabla{w^{\epsilon}}\phi_{-}|^p \\ &+&\sum_{k=1}^{m}C_{m}^{k}|\nabla{\phi}+\nabla{\phi_{-}}{w^{\epsilon}}|^{p-k} |\nabla{w^{\epsilon}}\phi_{-}|^{k} \\
&+ & \sum_{k=1}^{m}C_{m}^{k}|\nabla{\phi}+\nabla{\phi_{-}}{w^{\epsilon}}|^{k} |\nabla{w^{\epsilon}}\phi_{-}|^{p-k}.
\end{eqnarray*}

If we  use the same argument as above and we can get the same conclusion for $p$ is not an integer.

Thus for any $p$ : $1<p\leq n$, we have
$$\limsup_{\epsilon\rightarrow0}\int_{D}|\nabla{\phi}+\nabla{\phi_{-}}{w^{\epsilon}}+\nabla{w^{\epsilon}}\phi_{-}|^p dx \leq \int_{D}|\nabla \phi|^p dx + \int_{D}\alpha_{0}{\phi_{-}}^{p}dx$$

Hence
\begin{eqnarray*}
{\mathcal{F}}_{0}(\phi)&\geq&\limsup_{\epsilon\rightarrow0}{\mathcal{F}}(\phi+{\phi}_{-}w^{\epsilon})\\
&\geq&\liminf_{\epsilon\rightarrow0}{\mathcal{F}}({u^{\epsilon}})
\end{eqnarray*}

By  Key Lemma II, we have
$${\mathcal{F}}_{0}(\phi)\geq{\mathcal{F}}_{0}(u^{0}).  $$

And the set $\{\phi \in C^{1}_{0}(D): \phi_{-}\in C^{1}_{0}(D)\}$ is dense in $W^{1,p}_{0}(D)$, 
then $u^0$ is the solution of $$\min\{ \int_{D}\frac{1}{p}|\nabla{v}|^{p}+\frac{1}{p}\alpha_{0}v_{-}^{p}-fvdx : \  \forall \  v\in W_{0}^{1,p}(D)\}  \qquad \qquad \Box $$

{\bf Proof of Corollary 1.2 :} Let $\phi\in C^{1}_{0}(D) $
such that $\phi_{-}\in C^{1}_{0}(D)$ and $\phi\geq \psi$ a. e. in $D$. Then 
$$\phi_{+} + (w^{\epsilon}_{+}-1)\phi_{-} \ \geq \ \psi^{\epsilon} \ \mbox{in} \ D$$

Obviously, for $w^{\epsilon}_{+}$, we have the following property:
$$ \limsup_{\epsilon\rightarrow0}\int_{D}|\nabla w^{\epsilon}_{+}|^p |\phi_{-}|^p dx \leq \alpha_{0}\int_{D}|\phi_{-}|^p dx$$

Thus from the proof for Theorem 1.1 (Part II), we have that
\begin{eqnarray*}
&&\limsup_{\epsilon\rightarrow0}\int_{D}\frac{1}{p}|\nabla \{\phi_{+} + (w^{\epsilon}_{+}-1)\phi_{-} \}|^p dx\\
&=&\limsup_{\epsilon\rightarrow0}\int_{D}\frac{1}{p}|\nabla \phi_{+} +\nabla w^{\epsilon}_{+} \phi_{-} + (w^{\epsilon}_{+}-1)\nabla \phi_{-}|^pdx\\
&\leq&\int_{D}\frac{1}{p}|\nabla \phi|^p dx + \frac{\alpha_{0}}{p}\int_{D}(\phi_{-})^p dx
\end{eqnarray*}

And since $\phi_{+} + (w^{\epsilon}_{+}-1)\phi_{-} \ \geq \ \psi^{\epsilon} \ \mbox{in} \ D$, then
\begin{eqnarray*}
\int_{D}\frac{1}{p}|\nabla \{\phi_{+} + (w^{\epsilon}_{+}-1)\phi_{-} \}|^p dx&\geq& \int_{D}f\cdot(\phi_{+} + (w^{\epsilon}_{+}-1)\phi_{-})dx\\ 
&+&\int_{D}\frac{1}{p}|\nabla h^{\epsilon}|^p -f h^{\epsilon}dx 
\end{eqnarray*}

And $w^{\epsilon}_{+}$ converges to $0$ weakly in $W^{1,p}_{0}(D)$ and $h^{\epsilon}$ converges to $h_{0}$ weakly in $W^{1,p}_{0}(D)$, hence
\begin{eqnarray*}
 &&\liminf_{\epsilon\rightarrow0} \int_{D}\frac{1}{p}|\nabla h^{\epsilon}|^p -fh^{\epsilon}dx\\
&\leq& \int_{D}\frac{1}{p}|\nabla \phi|^p dx+ \frac{\alpha_{0}}{p}\int_{D}(\phi_{-})^p dx -f\phi dx
\end{eqnarray*}

By Lemma II, we have that
\begin{eqnarray*}
 &&\int_{D}\frac{1}{p}|\nabla h_{0}|^pdx +\int_{D}\frac{1}{p}\alpha_{0}(h_{0})_{-}^pdx -fh_{0} dx\\
&\leq& \int_{D}\frac{1}{p}|\nabla \phi|^p dx+ \frac{\alpha_{0}}{p}\int_{D}(\phi_{-})^p dx -\int_{D}f\phi dx
\end{eqnarray*}
And $\{\phi\in C^{1}_{0}(D): \ \phi_{-}\in C^{1}_{0}(D) \ \mbox{and}\ \phi\geq \psi \ \mbox{a.e.} \ \mbox{in} \ D\}$ is dense in $\{v\in W^{1,p}_{0}(D): \ v\geq \psi \ \mbox{a.e.} \ \mbox{in} \ D\}$,
therefore $h_0$ is the solution to the following variational problem:
$$\min\{\int_{D}\frac{1}{p}|\nabla{v}|^{p}+\frac{1}{p}\alpha_{0}v_{-}^{p}-fvdx: \ v\in W_{0}^{1,p}(D) \ \mbox{and} \ v\geq\psi \
 \mbox{a.}\ \mbox{e.}\ \mbox{in} \  D\}  .  \qquad \qquad   \Box$$

\section{ Proof of the Key Lemmas I} 
In order to prove Key Lemma I, we need to follow the the steps as follows:\\

{\bf Step 1:} Find the crital value $\alpha_{0}$.(Here we will use the method from Caffarelli and Mellet [CM].)

{\bf Step 2:} Show that $\{w^{\epsilon}\}$ is bounded in $W^{1,p}(D)$, where the correctors $\{w^{\epsilon}\}$ is defined by 
$$w_{\epsilon}(x,\omega)=\inf\{v(x): \ {\triangle}_{p}v \leq \alpha_{0} \ \mbox{in} \ D_{\epsilon}, 
\ v\geq1 \ \mbox{on} \ T_{\epsilon} \ \mbox{and} \ v=0\ \mbox{on} \ \partial{D}\}$$

{\bf Step 3:} Show that $w^{\epsilon}\longrightarrow0$ in $L^p(D)$ as $\epsilon\rightarrow0$.

{\bf Step 4:} Prove the property (a)-(c) of Key Lemma I,  i. e. 

{\bf(a)} for  any $\phi \in \mathcal{D}(D)$ and $0<p'<p$,
$$\lim_{\epsilon\rightarrow0}\int_{D}|\nabla w^{\epsilon}|^{p'}\phi dx =0$$

{\bf(b)} for  any $\phi \in \mathcal{D}(D)$,
$$\displaystyle\lim_{\epsilon\rightarrow0}\int_{D}|\nabla{w^{\epsilon}}|^{p}\phi dx=\int_{D}\alpha_{0}\phi dx .$$

{\bf(c)} for any sequence $\{v^{\epsilon}\}\subset W^{1,p}_{0}(D)$ 
with the property : $v^\epsilon \rightharpoonup v $ in $W^{1,p}_{0}(D)$ and $v^\epsilon =0$ on $T_{\epsilon}$ and 
any $\phi \in \mathcal{D}(D)$, we have that $$\lim_{\epsilon \rightarrow 0}\int_{D}|\nabla w^{\epsilon}|^{p-2}{\nabla w^{\epsilon}}\cdot {\nabla v^{\epsilon}}\phi dx = -\alpha_{0}\int_{D}v\phi dx.$$

\subsection{Find the crital value $\alpha_{0}$ }
Here we will use the method from Caffarelli and Mellet [CM]. we will consider the following obstacle problem:
for every open set $A\subset{\bf R}^n$ \ , $\alpha\in {\bf R}$
, $\forall x\in A$ and $\omega\in \Omega$ , we define
 $${v}^{\epsilon}_{\alpha,A}(x,\omega)=\inf\{v(x): \ {\triangle}_{p}{v}(\cdot)\leq \ \alpha-{\sum_{k\in {\bf Z}^{n}\cap {\epsilon}^{-1}A}\gamma(k,\omega)\epsilon^{n}\delta(\cdot  -\epsilon k)} \ \mbox{in} \ A,\  v\geq0\  \mbox{in} \ A,\  v=0 \  \mbox{on} \
 \partial{A}\}$$
where $1<p\leq n$ and 
 ${\triangle}_{p}{w}={\nabla}\cdot  ({\mid \nabla w \mid}^{p-2}{\nabla w})$ is the \ p-Laplacian \ operator
And we set
$${m}^{\epsilon}_{\alpha}(A, \omega)=|\{x\in A: \ {v}^{\epsilon}_{\alpha,A}=0\}|.$$

 From [CM], we can find that for any given $\epsilon>0$, \ 
 ${m}^{\epsilon}_{\alpha}(\cdot, \omega)$ is subadditive for each $\omega\in \Omega$ and the process ${m}^{\epsilon}_{\alpha}(A,\omega)$ is stationary ergodic.

Hence by [CM] and [DM], for any real number $\alpha$, there is a constant $l(\alpha)\geq0$ such that
$$\displaystyle\lim_{\epsilon\rightarrow 0}\frac{m^{\epsilon}_{\alpha}(B_{1}(x_{0}), \omega)}{|B_{1}(x_{0})|}=l(\alpha)\ ,$$
i.e.
$$\displaystyle\lim_{\epsilon\rightarrow 0}\frac{|\{x\in B_{1}(x_{0}):{v}^{\epsilon}_{\alpha,\ B_{1}(x_{0})}=0\}|}{|B_{1}(x_{0})|}=l(\alpha),$$
for any $B_{1}(x_{0})\subset {\bf R}^n$.

About the function $l(\alpha)$ , we have the following :

{\bf Propostion A.} 

 (i) \ \textit{ $l(\alpha)$ a nondecreasing function of $\alpha$ ;}

  (ii) \ \textit{$l(\alpha)=0 $ for $\alpha<0$;}

 (iii)\  \textit{ $l(\alpha)>0 $ for $\alpha$ is large enough.}

{\bf Proof:}\quad (i) \  For its monotonicity, we consider two parameters $\alpha\leq\alpha'$ and we will compare $l(\alpha)$ and $l(\alpha')$. By comparison principle, for any $A\subseteq {\bf R}^n$, $$v^{\epsilon}_{\alpha', A}(x, \omega)\leq v^{\epsilon}_{\alpha, A}(x, \omega),\ a.\ e.\ x\in A.$$ Hence $$\{x \in A: \ v^{\epsilon}_{\alpha, A}(x, \omega)=0\}\subseteq\{x \in A: \ v^{\epsilon}_{\alpha', A}(x, \omega)=0\},$$
 which implies $l(\alpha)\leq l(\alpha')$ for $\alpha\leq\alpha'$. \\
 \quad \ (ii) \ If $\alpha<0$, we let $\displaystyle \beta=|\alpha|^{(2-p)/(p-1)}\alpha$
then $${\triangle}_{p}\{\frac{\beta}{c(n,p)}|x-x_{0}|^{\frac{p}{p-1}}-\frac{\beta}{c(n,p)}\}= |\beta|^{p-2}\beta =\alpha$$ 
and $\frac{\beta}{c(n,p)}|x-x_{0}|^{\frac{p}{p-1}}-\frac{\beta}{c(n,p)}$ is positive in $B_{1}(x_{0})$ and vanishes on $\partial(B_{1}(x_{0}))$. Then we deduce that:
 $$v^{\epsilon}_{\alpha,  B_{1}}\geq \frac{\beta}{c(n,p)}|x-x_{0}|^{\frac{p}{p-1}}-\frac{\beta}{c(n,p)}>0 \ \mbox{in} \ B_{1}(x_{0}).$$

Hence $$m^{\epsilon}_{\alpha}(B_{1}(x_{0}), \omega)=0. $$
 Therefore $l(\alpha)=0$.\\
 \quad \  (iii)\ Let $a=a(k,\omega)=\sqrt[n]{\frac{nc\gamma(k,\omega)}{\alpha}}$, where the constant $c$ depends on $p$ and $n$.
 More precisely, it should be determined by the following: for $1<p<n$, 
$\displaystyle{{\triangle}_{p}}\{c^{\frac{1}{p-1}}|x|^{\frac{p-n}{p-1}}\}=-\delta(x)$; and for $p=n$, 
$\displaystyle{{\triangle}_{p}}\{c^{\frac{1}{n-1}}\log{\frac{1}{|x|}}\}=-\delta(x)$.

We  define the function $g^{\epsilon}_{\alpha, k}(x, \omega)$ for any $\alpha\in {\bf R}^n$ as following:
 \[ g^{\epsilon}_{\alpha, k}(x, \omega)=\left \{\begin{array}{ll}
\displaystyle{\int_{r}^{a\epsilon}(c\gamma(k,\omega){\epsilon}^{n}{s}^{1-n}-\frac{\alpha}{n}s)^{\frac{1}{p-1}}ds},& \ \mbox{if} \ 0\leq r=|x-\epsilon k|\leq a\epsilon;\\
0, & \ \mbox{if} \ x\in B_{1}\setminus B_{a\epsilon}(\epsilon k).
\end{array}
\right. \]
Obviously if the parameter $\alpha$ is large enough, then $\frac{1}{2}\geq \sqrt[n]{\frac{nc\gamma(k,\omega)}{\alpha}}$, which implies the function $g^{\epsilon}_{\alpha, k}(x, \omega)$ is only concentrated on the cell ball $B_{\frac{\epsilon}{2}}(\epsilon k)$ for $\alpha$ very large.
From the definition of $g^{\epsilon}_{\alpha, k}(x, \omega)$, we know that (for $\alpha$ is large)
$$\displaystyle{{\triangle}_{p}} \ g^{\epsilon}_{\alpha, k}(x, \omega)\leq \alpha -\gamma(k,\omega)\epsilon^{n}\delta(x-\epsilon k) \ \mbox{in} \ B_{1},$$
and $g^{\epsilon}_{\alpha, k}(x, \omega)=0$\ if \ $x\in B_{1}\setminus B_{a\epsilon}(\epsilon k)$.

Now we consider the sum of all $g^{\epsilon}_{\alpha, k}$:
$$\displaystyle\sum_{k\in \ {\epsilon}^{-1}B_{1}\cap {\bf Z}^n }g^{\epsilon}_{\alpha, k}  $$

By the definition, we know that for any two different $k, k'\in {\epsilon}^{-1}B_{1}\cap {\bf Z}^n$, $g^{\epsilon}_{\alpha, k}$ and $g^{\epsilon}_{\alpha, k'}$ have disjoint support. hence if we let $\displaystyle g_{\alpha}^{\epsilon}=\sum_{k\in \ {{\epsilon}^{-1}B_{1}}\cap {\bf Z}^n }g^{\epsilon}_{\alpha, k}$, then
$${{\triangle}_{p}} \ g^{\epsilon}_{\alpha}(x, \omega)\leq \alpha-\sum_{k\in \ {\epsilon}^{-1}B_{1}\cap {\bf Z}^n }\gamma(k,\omega)\delta(x-\epsilon k)$$

And  $g^{\epsilon}_{\alpha}(x, \omega)\geq 0 $ for $x\in B_{1}$\  and  $g^{\epsilon}_{\alpha}(x, \omega)=0$ on $\partial{B_{1}}$.

Therefore $$0\leq v^{\epsilon}_{\alpha, B_{1}}(x, \omega)\leq g^{\epsilon}_{\alpha}(x, \omega),\ \mbox{for} \ x\in B_{1}.$$

Thus $$\bigcup_{k\in \ {\epsilon}^{-1}B_{1}\cap {\bf Z}^n}(B_{1}\setminus B_{a\epsilon}(\epsilon k))\subset \{x\in B_{1}: \ v^{\epsilon}_{\alpha, B_{1}}=0\},$$
which implies $$m^{\epsilon}_{\alpha}(B_{1}, \omega)\geq \omega_{n}-C \epsilon^{-n}(a\epsilon)^{n}=1-Ca^{n},$$
where $\omega_{n}$ is the volume of the unit ball $B_{1}$.

And when $\alpha$ is large enough, then $a$ will be small enough such that $\omega_{n}-Ca^{n}\geq \frac{1}{2}\omega_{n}$. Therefore
$$m^{\epsilon}_{\alpha}(B_{1}, \omega)\geq \frac{1}{2}\omega_{n}>0,\ \mbox{if} \ \alpha \ \mbox{is} \ \mbox{large} \ \mbox{enough}$$

Then $l(\alpha)>0$ if $\alpha$ is large enough.  \qquad \qquad \qquad \qquad \qquad $\Box$

Next We choose the critical value $\alpha_{0}$ by the following way:
$$\alpha_{0}=\sup\{\alpha: l(\alpha)=0\}.$$
Then by {\bf Proposition A}, $\alpha_{0}$ is finite and nonnegative.    

In the following, we will define $w^{\epsilon}$ as follows:
$$w^{\epsilon}(x,\omega)=\inf\{v(x): \ {\triangle}_{p}v \leq \alpha_{0} \ \mbox{in} \ D_{\epsilon}, 
\ v\geq1 \ \mbox{on} \ T_{\epsilon} \ \mbox{and} \ v=0\ \mbox{on} \ \partial{D}\}$$
Therefore $w^{\epsilon}$ satisfies the following conditions:
\[ \left\{\begin{array}{lll}
\displaystyle{\triangle}_{p} w^{\epsilon}(x,\omega)&= \alpha_{0} \; \quad \mbox{for}\; \quad x\in D_{\epsilon}(\omega) \nonumber\\

w^{\epsilon}(x,\omega)&=1  \; \quad \mbox{for } \; \quad x \in T_{\epsilon}(\omega), \nonumber\\

w^{\epsilon}(x,\omega)&= 0 \; \quad \mbox{for} \; \quad x\in
{\partial D}\setminus T_{\epsilon}(\omega). \nonumber
\end{array}
\right.\] 

\subsection{ $W^{1,p}$ Boundedness of $\{w^\epsilon\}$}

To show that $\{w^{\epsilon}\}$ is uniformly bounded in $W^{1,p}(D)$,
 we split the proof into two parts: $\{w^{\epsilon}\}$ is uniformly bounded in $L^{p}(D)$ and
 $\{\nabla w^{\epsilon}\}$ is also uniformly bounded in $L^{p}(D)$.

To prove the first part: $\{w^{\epsilon}\}$ is uniformly bounded in $L^{p}(D)$, we need to introduce an auxiliary function $v(x)$: 
let $v$ be the solution to the following problem:
\[ \left\{\begin{array}{lll}
 {\triangle}_{p} v &= \alpha_{0} \ \mbox{in} \ D \\\
 v&=  0   \ \mbox{on}  \ {\partial D}    
\end{array}
\right. \]

By Comparison Principle for almost surely $\omega\in\Omega$, 
$$ v(x) \leq w^{\epsilon}(x, \omega)\leq 1 \ \mbox{for} \ \mbox{a. e.} \ x\in  {D}$$
Hence 
$$\int_{D}|w^{\epsilon}|^{p}dx \leq C$$
which implies that $\{w^{\epsilon}\}$ is uniformly bounded in $L^{p}(D)$

To show that $\{\nabla w^{\epsilon}\}$ is also uniformly bounded in $L^{p}(D)$, we define the function $h^{\epsilon}(x, \omega)$  as following: if $1<p<n$, then
\[h^{\epsilon}(x, \omega)=\left\{\begin{array}{lll}
1,  \qquad \qquad \qquad  \qquad \qquad \  x\in T_{\epsilon}\\
 \displaystyle \frac{(\frac{\epsilon}{2})^{\frac{p-n}{p-1}}-|x-\epsilon k|^{\frac{p-n}{p-1}}}{(\frac{\epsilon}{2})^{\frac{p-n}{p-1}}-(a^{\epsilon})^{\frac{p-n}{p-1}}}, \ \ x\in B_{\frac{\epsilon}{2}}\setminus B_{a^{\epsilon}}(\epsilon k)\\\ \\
0, \ \ \ \qquad \qquad \qquad \qquad \qquad \mbox{otherwise}
\end{array}
\right. \]
and if $p=n$, then
\[ h^{\epsilon}(x, \omega)=\left\{\begin{array}{lll}
1, \ \ \qquad \qquad \qquad \qquad  x\in T_{\epsilon}\\
 \displaystyle \frac{\log|x-\epsilon k| -\log{\frac{\epsilon}{2}}}{\log{a^{\epsilon}}-\log{\frac{\epsilon}{2}}}, \ \ x\in B_{\frac{\epsilon}{2}}\setminus B_{a^{\epsilon}}(\epsilon k)\\\ \\
0, \ \ \ \qquad \qquad \qquad \qquad \mbox{otherwise}
\end{array}
\right. \]

Obviously, $w^\epsilon -h^\epsilon =0$ on $T_{\epsilon}$ and $\partial{D}$. Hence
\begin{eqnarray*}
\alpha_{0}\int_{D}(h^\epsilon -w^\epsilon)dx &=& \int_{D}{\triangle}_{p} w^{\epsilon}(h^\epsilon -w^\epsilon)dx\\
&=&\int_{D}|\nabla w^{\epsilon}|^p dx -\int_{D}|\nabla w^{\epsilon}|^{p-2}\nabla w^{\epsilon}\cdot {\nabla h^{\epsilon}} dx
\end{eqnarray*}

Then by Holder inequality and Young's Inequality
\begin{eqnarray*}
\int_{D}|\nabla w^{\epsilon}|^p dx &\leq& (\int_{D}|\nabla w^{\epsilon}|^{p})^{\frac{p-1}{p}}\cdot (\int_{D}|\nabla h^{\epsilon}|^p dx)^{\frac{1}{p}}\\
&+& \alpha_{0}\int_{D}(h^\epsilon -w^\epsilon)dx\\
&\leq& \frac{p-1}{p}\int_{D}|\nabla w^{\epsilon}|^p dx + \frac{1}{p}\int_{D}|\nabla h^{\epsilon}|^p dx + C \alpha_{0}
\end{eqnarray*}

Thus 
\begin{eqnarray*}
 \int_{D}|\nabla w^{\epsilon}|^p dx &\leq& \int_{D}|\nabla h^{\epsilon}|^p dx + p \alpha_{0}\int_{D}|h^\epsilon -w^\epsilon|dx\\
&\leq& C 
\end{eqnarray*}
where C is a universal constant depending on n and $\alpha_{0}$.

Therefore  $$\{w^{\epsilon}\} \ \mbox{ is uniformly bounded in}  W^{1,p}(D). \qquad \qquad \qquad \qquad \Box$$

\subsection{$w^{\epsilon}\longrightarrow0$ in $L^{p}(D)$ as $\epsilon\rightarrow0$}
To prove this fact: $\displaystyle\lim_{\epsilon\rightarrow0}\int_{D}|w^{\epsilon}|^p dx =0,$
we need to compare $w^{\epsilon}$ with ${v}^{\epsilon}_{\alpha_{0}, D}$. Roughly speaking,  we will show that, 
near the singular points (the holes), their liming behaviour should be very close.

First of all , we will consider some asymptotic properties of ${v}^{\epsilon}_{\alpha_{0}, D}$ as $\epsilon \rightarrow 0$.
Since $D$ is bounded, then without loss of generality we can assume that $D\subset B_{1}$.
In the following we will use ${v}^{\epsilon}_{0}$  and $\overline{v}^{\epsilon}_{0}$ to  denote ${v}^{\epsilon}_{\alpha_{0}, B_{1}}$ and $\displaystyle \min({v}^{\epsilon}_{\alpha_{0}, B_{1}},\ 1)$. 
And about ${v}^{\epsilon}_{0}$, we have the following facts:

{\bf Proposition B.} (i) \textit{  we have that ${v}^{\epsilon}_{\alpha_{0}, D}(x,\omega)\geq h^{\epsilon}_{k}(x, \omega)-o(1) \ \ 
for \  a.\ e.\ x\in B_{\frac{\epsilon}{2}}(\epsilon k)  \ and \ a.\ s.\ \omega\in\Omega$,
where 
\[ h^{\epsilon}_{k}(x, \omega)=\left
\{\begin{array}{ll}
 & c^{\frac{1}{p-1}} \gamma(k, \omega)^{\frac{1}{p-1}}\epsilon^{\frac{n}{p-1}}|x-\epsilon k|^{\frac{p-n}{p-1}},   \qquad  \mbox{if} \qquad 1<p<n, \nonumber\\\
 & -c^{\frac{1}{n-1}}\gamma(k, \omega)^{\frac{1}{n-1}} \epsilon^{\frac{n}{n-1}}\log{|x-\epsilon k|} ,     \qquad \mbox{if} \qquad p=n. \nonumber
\end{array}
\right. \]
and the constant $c$ is the same constant as Propostion A (iii).}

(ii) \textit{ For any $\delta >0$, $\overline{v}^{\epsilon}_{\delta}(x,\omega) \ converges \ to \ 0 \ in \ L^p(B_{1})  \ as \ \epsilon \ goes \ to \ 0\ \ 
for \ a.\ s.\ \omega\in\Omega$, where $\overline{v}^{\epsilon}_{\delta}$ is defined as follows:
\begin{eqnarray*} v^{\epsilon}_{\delta}(x,\omega)&=&\inf\{v(x): {\triangle}_{p}{v}\leq \alpha_{0} +  \delta\ - \sum_{k\in {\bf Z}^{n}\cap {\epsilon}^{-1}B_{1}}\gamma(k,\omega)\epsilon^{n}\delta(\cdot -\epsilon k) \\ 
&\mbox{in}& \quad B_{1}, \ v\geq0 \quad \mbox{on} \quad B_{1}, \ v=0\ \mbox{on} \quad \partial{B_{1}}\}
\end{eqnarray*}
 and let $\overline{v}^{\epsilon}_{\delta}=\min(v^{\epsilon}_{\delta},\  1)$. 
Hence
$$\overline{v}^{\epsilon}_{\alpha_{0} +\delta, D}=\min({v}^{\epsilon}_{\alpha_{0} +\delta, D}, \ 1)  \ converges \ to \ 0 \ in \ L^p(D) $$ }

{\bf Proof:}  ( i) Let $\displaystyle b(k,\omega)=\sqrt[n]{\frac{nc\gamma(k,\omega)}{\alpha_{0}}}$. 
Then we define the function $h^{\epsilon}_{\alpha_{0}, k}(x, \omega)$ as follows: if $b\geq \frac{1}{2}$, then
\[ h^{\epsilon}_{\alpha_{0}, k}(x, \omega)= \left\{\begin{array}{lll}
\displaystyle{\int_{r}^{\frac{\epsilon}{2}}(c\gamma(k,\omega){\epsilon}^{n}{s}^{1-n}-\frac{\alpha_{0}}{n}s)^{\frac{1}{p-1}}ds},& 0\leq r\leq\frac{1}{2}\epsilon\\
0 , &  r \geq \frac{\epsilon}{2}
\end{array}
\right. \]
and if $b\leq \frac{1}{2}$, then
\[ h^{\epsilon}_{\alpha_{0}, k}(x, \omega)= \left\{\begin{array}{lll}
\displaystyle{\int_{r}^{b{\epsilon}}(c\gamma(k,\omega){\epsilon}^{n}{s}^{1-n}-\frac{\alpha_{0}}{n}s)^{\frac{1}{p-1}}ds},& 0\leq r\leq b\epsilon\\
0 , &  r \geq b{\epsilon}
\end{array}
\right. \]
where $r=|x-\epsilon k|$ and $x\in B_{\frac{1}{2}\epsilon}(\epsilon k)$.

If $b\geq \frac{1}{2}$ , we have that $\forall \  x \in B_{b_{\epsilon}}(\epsilon k) \ \mbox{and} \ \ \omega\in\Omega$,
$${\triangle}_{p} h^{\epsilon}_{\alpha_{0}, k}(x, \omega)=\alpha_{0}-\gamma(k,\omega){\epsilon}^{n}\delta(x-\epsilon k)  \ \mbox{in} \ B_{\frac{\epsilon}{2}}(\epsilon k)$$
and $h^{\epsilon}_{\alpha, k}(x, \omega)=0$ if $|x-\epsilon k|=\frac{\epsilon}{2}$.
Hence
$${\triangle}_{p} h^{\epsilon}_{\alpha_{0}, k}(x, \omega)\geq{\triangle}_{p} {v}^{\epsilon}_{\alpha_{0},D}(x,\omega),\  \forall \ x \in B_{\frac{\epsilon}{2}}(\epsilon k)$$
and $ h^{\epsilon}_{\alpha_{0}, k}(x, \omega)=0\leq {v}^{\epsilon}_{0}(x,\omega)$ when $|x-\epsilon k|=\frac{\epsilon}{2}$.
By comparison principle (see [MZ]), we have $$h^{\epsilon}_{\alpha_{0}, k}(x, \omega)\leq {v}^{\epsilon}_{\alpha_{0},D}(x,\omega),\ \mbox {a.\ e.} \ x \in B_{\frac{\epsilon}{2}}(\epsilon k) \ \mbox{and} \  \mbox {a.\ s.} \ \omega\in\Omega$$
And for the case $b\leq \frac{1}{2}$, the proof is similar as above.

Thus
$$h^{\epsilon}_{\alpha_{0}, k}(x, \omega)\leq {v}^{\epsilon}_{\alpha_{0},D}(x,\omega),\ \mbox{a.\ e.} \ x \in B_{\frac{\epsilon}{2}}(\epsilon k) \ \mbox{and} \ \mbox {a.\ s.} \ \omega\in\Omega$$
And by direct simple computation, we know that 
$$h^{\epsilon}_{\alpha_{0}, k}(x, \omega)\geq h^{\epsilon}_{k}(x, \omega)-o(1), \ \mbox{a. e.}\  x\in B_{\frac{\epsilon}{2}}(\epsilon k) \ \mbox{and} \ \mbox {a.\ s.}\ \omega\in\Omega$$
Therefore
 $$\displaystyle{v}^{\epsilon}_{\alpha_{0},D}(x,\omega)\geq h^{\epsilon}_{k}(x, \omega)-o(1), \ \mbox{a. e.}\  x\in B_{\frac{\epsilon}{2}}(\epsilon k) \ \mbox{and} \ \mbox {a.\ s.} \ \omega\in\Omega$$
which concludes Proposition B (i).

(ii) From the definition of $\{\overline{v}^{\epsilon}_{\delta}\}$, we know that for a. e. $\omega \in \Omega$
$$-\ \sum_{k\in {\bf Z}^{n}\cap {\epsilon}^{-1}B_{1}}\gamma(k,\omega)\epsilon^{n}\delta(\cdot -\epsilon k) \leq{\triangle}_{p}{{v}^{\epsilon}_{\delta}}\leq \alpha_{0}\ + \ \delta \qquad \mbox{in} \quad B_{1}$$

Hence $$\langle {\triangle}_{p}{{v}^{\epsilon}_{\delta}}, \ \overline{v}^{\epsilon}_{\delta} \rangle \geq -\ \sum_{k\in {\bf Z}^{n}\cap {\epsilon}^{-1}B_{1}}\gamma(k,\omega)\epsilon^{n}\overline{v}^{\epsilon}_{\delta}(\epsilon k)$$

By Proposition B (i), we have that 
$${v}^{\epsilon}_{\delta}\geq h^{\epsilon}_{k}(x, \omega)-o(1), \ \mbox{a. e.}\  x\in B_{\frac{\epsilon}{2}}(\epsilon k) \ \mbox{and} \ \mbox {a.\ s.} \ \omega\in\Omega$$
which concludes that
$\displaystyle\overline{v}^{\epsilon}_{\delta}(\epsilon k)=1$

Thus from integration by parts
$$\int_{D}|\nabla \overline{v}^{\epsilon}_{\delta} |^p dx \leq C$$
where $C$ is a universal constant.
Therefore $\{\overline{v}^{\epsilon}_{\delta}\}$  are bounded in $W^{1,p}(B_{1})$.

From [CM] and [CSW], we have that for a. s. $\omega\in\Omega$.
$$\lim_{\epsilon\rightarrow0}\frac{|\{\overline{v}^{\epsilon}_{\delta}=0\}\cap B_{r}(x_{0})|}{|B_{r}(x_{0})|}=l(\alpha_{0}+\delta)>0, \ \mbox{for} \ \mbox{any} \
B_{r}(x_{0})\subseteq B_{1}.$$
By the Poincar\'{e}-Sobolev inequality (see Lemma 4.8 in [HL]), there exists a constant $C=C(\alpha_{0}+\delta,n)$ such that
$$\int_{B_{r}(x_{0})}|\overline{v}^{\epsilon}_{\delta}|^p dx \leq C r^{p}\int_{B_{r}(x_{0})}|\nabla \overline{v}^{\epsilon}_{\delta}|^p dx$$
for any $B_{r}(x_{0})\subseteq B_{1}$

And $\{\overline{v}^{\epsilon}_{\delta}\}$  are bounded in $W^{1,p}(B_{1})$, hence
$$\int_{B_{1}}|\overline{v}^{\epsilon}_{\delta}|^p dx \leq C r^{p}$$
which implies that
$$\lim_{\epsilon\rightarrow0}\int_{B_{1}}|\overline{v}^{\epsilon}_{\delta}|^p dx=0$$
 And $\overline{v}^{\epsilon}_{\delta}(x, \omega)\geq \overline{v}^{\epsilon}_{\alpha_{0}+\delta, \ D}(x, \omega) \geq 0$ for a. e. $x\in D$ and a. s. $\omega\in\Omega$ , thus
$$\lim_{\epsilon\rightarrow0}\int_{D}|\overline{v}^{\epsilon}_{\alpha_{0}+\delta, D}|^p dx=0 \qquad \qquad \qquad \qquad \Box
$$

To finish the proof of {\bf Step 3 }, we need to  pass the limiting property of $\overline{v}^{\epsilon}_{\alpha_{0}+\delta, D}$ to $w^\epsilon$. 
First of all, we introduce a new auxillary function $w^{\epsilon}_{\delta}$ as follows: 
$$w^{\epsilon}_{\delta}(x,\omega)=\inf\{v(x): \ {\triangle}_{p}v \leq \alpha_{0} + \delta \ \mbox{in} \ D_{\epsilon}, 
\ v\geq1 \ \mbox{on} \ T_{\epsilon} \ \mbox{and} \ v=0\ \mbox{on} \ \partial{D}\}$$

Obviously, for a. s. $\omega \in \Omega$,  $w^{\epsilon}(x, \omega)\geq w^{\epsilon}_{\delta}(x, \omega)$ for a. e. $x\in D$  
and $\{w^{\epsilon}_{\delta}\}$ is also bounded in $W^{1,p}(D)$ by {\bf Step 2} 
(the argument for the $W^{1,p}$ boundedness of $w^\epsilon$).

More precisely, $\{w^{\epsilon}_{\delta}\}$ satisfies the following property:

{\bf Proposition C.} (i) \textit{ For any $\delta>0$, 
\[ \|w^{\epsilon}-w^{\epsilon}_{\delta}\|_{W^{1,p}(D)}\leq\left\{\begin{array}{lll}
C \delta^{1/(p-1)}, \ \ 2 \leq p\leq n \\
C \delta, \qquad \quad 1<p\leq 2
\end{array}
\right.\]
where C depends only on $p$, $n$ and $\alpha_{0}$}

(ii) \textit{$$\lim_{\epsilon\rightarrow0}\int_{D}{(w^{\epsilon}_{\delta})^p} dx =0$$}

{\bf Proof of Proposition C:} (i) Next we apply the well-known inequality
\[(|\xi|^{p-2}\xi-|\eta|^{p-2}\eta)\cdot(\xi-\eta)\geq \gamma\left\{\begin{array}{lll}
|\xi-\eta|^2(|\xi|+|\eta|)^{p-2}, \quad 1<p\leq 2 \\
|\xi-\eta|^p, \qquad \qquad  \qquad 2\leq p\leq n
\end{array}
\right.\]
for any nonzero $\xi,\ \eta \in {\bf R}^n$ and a constant $\gamma=\gamma(n,p)>0$.

If $2\leq p\leq n$, then we have the following:
\begin{eqnarray*}
\displaystyle&&\int_{D}\{{\triangle}_{p}  w^{\epsilon}_{\delta}(x,\omega)-{\triangle}_{p} w^{\epsilon}(x,\omega)\}\cdot\{w^{\epsilon}(x,\omega)-w^{\epsilon}_{\delta}(x,\omega)\}dx\\
&=&\int_{D}\{|\nabla w^{\epsilon}_{\delta}|^{p-2}\nabla w^{\epsilon}_{\delta}-|\nabla w^{\epsilon}|^{p-2}\nabla w^{\epsilon}\}\cdot \{\nabla w^{\epsilon}_{\delta}-\nabla w^{\epsilon}\}dx\\
&\geq& \gamma\int_{D}|\nabla w^{\epsilon}_{\delta}-\nabla w^{\epsilon}|^{p}dx.
\end{eqnarray*}

And by H\"{o}lder inequality and Poincar\'{e} inequality, we have
\begin{eqnarray*}
\displaystyle\int_{D}\delta \{w^{\epsilon}(x,\omega)-w^{\epsilon}_{\delta}(x,\omega)\}&dx&\leq C \delta \{\int_{D}|w^{\epsilon}(x,\omega)-w^{\epsilon}_{\delta}(x,\omega)|^p dx\}^{\frac{1}{p}}\\
&\leq& C\delta \{\int_{D}|\nabla w^{\epsilon}(x,\omega)-\nabla w^{\epsilon}_{\delta}(x,\omega)|^p dx\}^{\frac{1}{p}}
\end{eqnarray*}
Therefore $$\displaystyle\int_{D}|\nabla w^{\epsilon}_{\delta}-\nabla w^{\epsilon}|^{p}dx\leq C\delta^{\frac{p}{p-1}}$$
which implies that$$\|w^{\epsilon}_{\delta}-w^{\epsilon}\|_{W^{1,p}(D)}\leq C \delta^{\frac{1}{p-1}}$$
where $C$ depends only on $p, \ n$.

If $1<p\leq 2$, then by H\"{o}lder inequality
\begin{eqnarray*}
&&\int_{D}\{{\triangle}_{p}  w^{\epsilon}_{\delta}-{\triangle}_{p} w^{\epsilon}\}\cdot\{w^{\epsilon}-w^{\epsilon}_{\delta}\}dx\\
&=&\int_{D}\delta \{w^{\epsilon}(x,\omega)-w^{\epsilon}_{\delta}(x,\omega)\}dx\\
&\geq& \gamma \int_{D}|\nabla(w^{\epsilon}-w^{\epsilon}_{\delta})|^2 (|\nabla w^{\epsilon}|+|\nabla w^{\epsilon}_{\delta}|)^{p-2} dx\\
&\geq&\gamma (\int_{D}|\nabla(w^{\epsilon}-w^{\epsilon}_{\delta})|^p dx)^{2/p}\times (\int_{D} (|\nabla w^{\epsilon}|+|\nabla w^{\epsilon}_{\delta}|)^p dx)^{1-2/p}
\end{eqnarray*}
And by H\"{o}lder inequality and Poincar\'{e} inequality,
$$\int_{D}\delta \{w^{\epsilon}(x,\omega)-w^{\epsilon}_{\delta}(x,\omega)\}dx
\leq C\delta \{\int_{D}|\nabla w^{\epsilon}(x,\omega)-\nabla w^{\epsilon}_{\delta}(x,\omega)|^p dx\}^{1/p}$$

Then
\begin{eqnarray*}
&&\{\int_{D}|\nabla w^{\epsilon}(x,\omega)-\nabla w^{\epsilon}_{\delta}(x,\omega)|^p dx\}^{1/p}\\
&\leq& C\delta (\int_{D} (|\nabla w^{\epsilon}|+|\nabla w^{\epsilon}_{\delta}|)^p dx)^{2/p-1}\\
&\leq& C\delta(\int_{D}|\nabla w^{\epsilon}|^p dx +\int_{D}|\nabla w^{\epsilon}_{\delta}|^p dx)^{2/p-1}\\
&\leq&C\delta
\end{eqnarray*}
Therefore 
\[ \|w^{\epsilon}-w^{\epsilon}_{\delta}\|_{W^{1,p}(D)}\leq\left\{\begin{array}{lll}
C \delta^{1/(p-1)}, \ \ 2 \leq p\leq n \\
C \delta, \qquad \quad 1<p\leq 2
\end{array}
\right.\]
where C depends only on $p$, $n$ and $\alpha_{0}$.

(ii)
And by Proposition B (i) and comparison principle, 
we know that $$ 0\leq (w^{\epsilon}_{\delta})_{+} \leq \overline{v}_{\alpha_{0}+\delta, \ D}^{\epsilon} + o(1)$$
hence $$\lim_{\epsilon\rightarrow0}\int_{D}{(w^{\epsilon}_{\delta})^p_{+}} dx =0$$
Next we will consider $(w^{\epsilon}_{\delta})_{-}$. 
In $\displaystyle B_{\epsilon/2}(\epsilon k)$, we suppose that $$\sup_{B_{\epsilon/2}(\epsilon k)}(w^{\epsilon}_{\delta})_{-}>0$$

Since ${\triangle}_{p}  w^{\epsilon}_{\delta}=\alpha_{0}+\delta$ in $D_{\epsilon}$, then 
by [MZ] $ w^{\epsilon}_{\delta}$ is continuous in $D$ and so is $(w^{\epsilon}_{\delta})_{-}$.
Then if  we  apply Harnack inequality (see [MZ]) to $(w^{\epsilon}_{\delta})_{-}$, we will have that
 for a. s. $\omega\in \Omega$,
  $$\sup_{B_{\epsilon/2}(\epsilon k)}(w^{\epsilon}_{\delta})_{-} = o(1),\ \mbox{for}\ \epsilon \ \mbox{is \ small}$$
which implies that
$$\lim_{\epsilon\rightarrow0}\int_{D}{(w^{\epsilon}_{\delta})^p_{-}} dx =0$$
Thus
$$\lim_{\epsilon\rightarrow0}\int_{D}{|w^{\epsilon}_{\delta}|^p} dx =0$$
which concludes Proposition C.   \qquad   \qquad    \qquad \qquad \qquad \qquad \qquad $\Box$

Hence by Proposition C, we have that
$$\lim_{\epsilon\rightarrow0}\int_{D}{|w^{\epsilon}|^p} dx =0$$
Therefore we can select a subsequence from $\{w^{\epsilon}\}$ such that 
this subsequence (we still use $\{w^{\epsilon}\}$ to denote it) converges weakly to zero in $W^{1,p}(D)$.

\subsection{Property (a)-(c) of $w^\epsilon$}
  
{\bf (a)} \ Without loss of the generality, we assume that $\phi\in C^{1}_{0}(D)$  and $\phi \geq0$ on $D$. 
Let $\theta$ be an any small positive number ($0<\theta<1$).To prove  property (a), 
we need to prove the two facts:
 $$ \limsup_{\epsilon\rightarrow0}\int_{D\cap \{w^{\epsilon}\leq \theta\}}|\nabla w^{\epsilon}|^{p'} \phi dx \leq C(\alpha_{0}, \phi)\theta, $$
and 
$$ \limsup_{\epsilon\rightarrow0}\int_{D\cap \{w^{\epsilon}> \theta\}}|\nabla w^{\epsilon}|^{p'} \phi dx =0 $$ 
In fact we let $w^{\epsilon}_{\theta} =(\theta -w^\epsilon)_{+}$, 
then $w^{\epsilon}_{\theta} \in W^{1,p}_{0}(D)$ and $w^{\epsilon}_{\theta}$ converges to $\theta$ weakly in $W^{1,p}_{0}(D)$. And since $\theta<1$, then $w^{\epsilon}_{\theta}$ on the holes $T_\epsilon$.
Obviously $$\lim_{\epsilon\rightarrow0}\int_{D}|\nabla w^{\epsilon}|^{p-2}\nabla w^{\epsilon}\cdot \nabla (w^{\epsilon}_{\theta}\phi) dx =-\alpha_{0}\theta \int_{D}\phi dx$$
which implies
\begin{eqnarray*}
 \lim_{\epsilon\rightarrow0}\{\int_{D\cap \{w^\epsilon\leq \theta\}}|\nabla w^\epsilon|^p \phi dx& -& \int_{D\cap \{w^\epsilon\leq \theta\}}|\nabla w^{\epsilon}|^{p-2}\nabla w^{\epsilon}\cdot \nabla \phi \ w^{\epsilon}_{\theta} dx \}\\
&=& \alpha_{0}\theta\int_{D}\phi dx
\end{eqnarray*}

Since $w^{\epsilon}$ is bounded in $W^{1,p}$, then by H\"{o}lder inequality,
\begin{eqnarray*}
|\int_{D\cap \{w^{\epsilon}\leq \theta\}}|\nabla w^{\epsilon}|^{p-2}\nabla w^{\epsilon}\cdot \nabla \phi \ w^{\epsilon}_{\theta} dx| &\leq& C \{\int_{D\cap \{w^{\epsilon}\leq \theta\}}(w^{\epsilon}_{\theta})^p dx\}^{1/p}\\
&\leq& C \{\int_{D}(w^{\epsilon}_{\theta})^p dx\}^{1/p}
\end{eqnarray*}
And $w^{\epsilon}$ converges to $0$ in $W^{1,p}$ weakly, then 
$$\lim_{\epsilon\rightarrow0}\{\int_{D} (w^{\epsilon}_{\theta})^p dx\}^{\frac{1}{p}}=\theta$$
Thus
$$\limsup_{\epsilon\rightarrow0}\int_{D\cap \{w^\epsilon\leq \theta\}}|\nabla w^\epsilon|^p \phi dx \leq C\theta$$

Now let $0<p'<p$, then by H\"{o}lder inequality, we have that
$$\int_{D\cap \{w^{\epsilon}\leq \theta\}}|\nabla w^{\epsilon}|^{p'} \phi dx \leq\{\int_{D\cap \{w^{\epsilon}\leq \theta\}}|\nabla w^{\epsilon}|^p \phi dx \}^{\frac{p'}{p}}\cdot \{\int_{D}\phi dx\}^{\frac{p-p'}{p}}$$
Thus $$ \limsup_{\epsilon\rightarrow0}\int_{D\cap \{w^{\epsilon}\leq \theta\}}|\nabla w^{\epsilon}|^{p'} \phi dx \leq C(\alpha_{0}, \phi)\theta,      $$

And for the integral $\displaystyle\int_{D\cap \{w^{\epsilon}> \theta\}}|\nabla w^{\epsilon}|^{p'} \phi dx $, 
we still apply H\"{o}lder inequality, then 
$$\int_{D\cap \{w^{\epsilon}> \theta\}}|\nabla w^{\epsilon}|^{p'} \phi dx \leq\{\int_{D}|\nabla w^{\epsilon}|^p \phi dx \}^{\frac{p'}{p}}\cdot \{\int_{D\cap \{w^{\epsilon}>\theta\}}\phi dx\}^{\frac{p-p'}{p}}$$

And $w^\epsilon \rightharpoonup 0 $ weakly in $W^{1,p}_{0}(D)$, then $$\lim_{\epsilon\rightarrow0}\int_{D\cap \{w^{\epsilon}>\theta\}}\phi dx =0$$
Then
$$ \limsup_{\epsilon\rightarrow0}\int_{D\cap \{w^{\epsilon}> \theta\}}|\nabla w^{\epsilon}|^{p'} \phi dx =0      $$
Therefore
$$\limsup_{\epsilon\rightarrow0}\int_{D}|\nabla w^{\epsilon}|^{p'} \phi dx \leq C(\alpha_{0}, \phi)\theta$$

And $\theta $ is an arbitrary small positive number, so

$$\lim_{\epsilon\rightarrow0}\int_{D}|\nabla w^{\epsilon}|^{p'} \phi dx =0. \qquad \qquad \qquad \qquad \qquad\Box$$

{\bf(b)} By integration by parts, we have:
\begin{eqnarray*}
\displaystyle \alpha_{0}\int_{D}\phi(1-w^{\epsilon})dx&=&\int_{D}{\nabla}\cdot  ({\mid \nabla w^{\epsilon} \mid}^{p-2}{\nabla w^{\epsilon}})\phi(1-w^{\epsilon})dx\\
&=&\int_{D}{\nabla \phi}\cdot{\nabla w^{\epsilon}}|\nabla w^{\epsilon}|^{p-2}(w^{\epsilon}-1)dx\\
&+&\int_{D}\phi|\nabla w^{\epsilon}|^{p}dx
\end{eqnarray*}
Since $w^{\epsilon}$ goes to 0 weakly in $W^{1,p}(D)$, hence 
$$\lim_{\epsilon\rightarrow0}\alpha_{0}\int_{D}\phi(1-w^{\epsilon})dx=\alpha_{0}\int_{D}\phi dx$$
And $w^{\epsilon}$ converges to 0 strongly in $L^{p}(D)$ and $\nabla w^{\epsilon}$ is bounded in $L^{p}(D)$ , 
hence by H\"{o}lder inequality, we have that
$$\lim_{\epsilon\rightarrow0}\int_{D}{\nabla \phi}\cdot{\nabla w^{\epsilon}}|\nabla w^{\epsilon}|^{p-2}w^{\epsilon}dx =0$$

Finally, by (a), we know that 
$$\lim_{\epsilon\rightarrow0}\int_{D}{\nabla \phi}\cdot{\nabla w^{\epsilon}}|\nabla w^{\epsilon}|^{p-2}dx =0$$

Therefore
$$\lim_{\epsilon\rightarrow0}\int_{D}\phi |\nabla w^{\epsilon}|^{p}dx=\int_{D}\alpha_{0}\phi dx.  \qquad \qquad \qquad \qquad \Box$$
{\bf (c)} From integration by parts, we have that
\begin{eqnarray*}
 \int_{D}{\nabla}\cdot  (| \nabla w^{\epsilon} |^{p-2}{\nabla w^{\epsilon}})v^{\epsilon}\phi dx&=& -\int_{D}\phi| \nabla w^{\epsilon} |^{p-2}{\nabla w^{\epsilon}}\cdot \nabla{v^{\epsilon}}dx\\
&-&\int_{D}v^{\epsilon}|\nabla w^\epsilon|^{p-2}\nabla w^\epsilon \cdot \nabla \phi dx
\end{eqnarray*}

which concludes that
$$-\int_{D}\alpha_{0}v^\epsilon \phi = \int_{D}\phi| \nabla w^{\epsilon} |^{p-2}{\nabla w^{\epsilon}}\cdot \nabla{v^{\epsilon}}dx+\int_{D}v^{\epsilon}|\nabla w^\epsilon|^{p-2}\nabla w^\epsilon \cdot \nabla \phi dx$$
Since  $v^\epsilon $  is bounded in $W^{1,p}_{0}(D)$ ($1<p\leq n$), 
then by Sobolev imbedding theorem ([GT]) $v^\epsilon$ is bounded in $L^{q}$ for some $q>p$.
Hence by H\"{o}lder inequality,
$$|\int_{D}v^{\epsilon}|\nabla w^\epsilon|^{p-2}\nabla w^\epsilon \cdot \nabla \phi dx|\leq \{\int_{D}|v^{\epsilon}|^q\}^{\frac{1}{q}} \{\int_{D}|\nabla w^\epsilon|^{(p-1)q'}|\nabla \phi|^{q'}dx\}^{\frac{1}{q'}}$$
where $q'=\frac{q}{q-1}$. And $q>p$, then  $(p-1)q'< p$, which implies that (by (a))
$$\lim_{\epsilon\rightarrow0}\{\int_{D}|\nabla w^\epsilon|^{(p-1)q'}|\nabla \phi|^{q'} dx\}^{\frac{1}{q'}}=0$$
Hence 
$$\lim_{\epsilon\rightarrow0}\int_{D}v^{\epsilon}|\nabla w^\epsilon|^{p-2}\nabla w^\epsilon \cdot \nabla \phi dx=0$$
Therefore
$$\lim_{\epsilon \rightarrow 0}\int_{D}|\nabla w^{\epsilon}|^{p-2}{\nabla w^{\epsilon}}\cdot {\nabla v^{\epsilon}}\phi dx = -\alpha_{0}\int_{D}v\phi dx.  \qquad \qquad \qquad \Box$$

\section{Proof of Key Lemma II }

{\bf Proof of Key Lemma II:} \ Let us decompose $u^{\epsilon}=u^{\epsilon}_{+}-u^{\epsilon}_{-}$, where $u^{\epsilon}_{+} = \max\{u^\epsilon, 0\}$ and $u^{\epsilon}_{-} = \max\{-u^{\epsilon},0\}$. Since $u^\epsilon$ converges to $u_{0}$ weakly in $W^{1,p}$, then $\displaystyle u^{\epsilon}_{+}\rightharpoonup u^{0}_{+}$ in $W^{1,p}$  ( $\displaystyle u^{\epsilon}_{-}\rightharpoonup u^{0}_{-}$ in $W^{1,p}$, respectively). 
Obviously,  $\displaystyle\int_{D}|\nabla u^{\epsilon}|^p dx=\int_{D}|\nabla u^{\epsilon}_{+}|^pdx+\int_{D}|\nabla u^{\epsilon}_{-}|^pdx$ and $\displaystyle\int_{D}|\nabla u^{0}|^pdx=\int_{D}|\nabla u^{0}_{+}|^pdx+\int_{D}|\nabla u^{0}_{-}|^pdx$

For $u^{\epsilon}_{+}$, we apply the classical lower semicontinuity property:
$$\liminf_{\epsilon\rightarrow0}\int_{D}|\nabla u^{\epsilon}_{+} |^p \geq \int_{D}|\nabla u^{0}_{+}|^p$$

In order to prove Lemma II, we need to prove the following revised lower semicontinuity property:
$$ \liminf_{\epsilon\rightarrow0}\int_{D}|\nabla u^{\epsilon}_{-}|^p \geq \int_{D}|\nabla u^{0}_{-}|^pdx + \int_{D}\alpha_{0}(u^{0}_{-})^p dx$$

Let $\theta $ be an any (small) positive number and $\phi$ is a test function (which is in $C^{1}_{0}(D)$ ).

Firstly we will show that
$$
\liminf_{\epsilon\rightarrow0}\frac{1}{p}\int_{w^\epsilon \leq \theta}|\nabla u^{\epsilon}_{-} |^pdx \geq\int_{D} |\nabla \phi|^{p-2}\nabla \phi \cdot \nabla u_{0}^{-}dx
-\frac{p-1}{p}\int_{D}|\nabla \phi|^p dx \qquad (*)$$

In fact from Young's inequality, we have the following
$$\int_{w^\epsilon \leq \theta}|\nabla \phi|^{p-2}\nabla \phi \cdot \nabla u^{\epsilon}_{-} dx\leq\frac{1}{p}\int_{w^\epsilon \leq \theta}|\nabla u^{\epsilon}_{-} |^p dx+ \int_{w^\epsilon \leq \theta}\frac{p-1}{p}|\nabla \phi|^p dx.$$

Since $w^\epsilon$ converges to 0 weakly in $W^{1,p}(D)$, then $\displaystyle |\{w^\epsilon >\theta\}|\rightarrow0$
as $\epsilon $ goes to 0. Hence

$$\lim_{\epsilon\rightarrow0}\int_{w^\epsilon >\theta}|\nabla \phi|^pdx=0$$
which implies that ( by H\"{o}lder inequality )
$$\lim_{\epsilon\rightarrow0}\int_{w^\epsilon >\theta} |\nabla \phi|^{p-2}\nabla \phi \cdot \nabla u^{\epsilon}_{-}  dx=\lim_{\epsilon\rightarrow0}\int_{w^\epsilon >\theta}|\nabla \phi|^p dx =0$$
Since $u^{\epsilon}_{-} $ converges to $u^{0}_{-}$ weakly in $W^{1,p}(D)$, then we have the estimate $(*)$:
$$\liminf_{\epsilon\rightarrow0}\frac{1}{p}\int_{w^\epsilon \leq \theta}|\nabla u^{\epsilon}_{-} |^pdx \geq \int_{D} |\nabla \phi|^{p-2}\nabla \phi \cdot \nabla u_{0}^{-}dx -\frac{p-1}{p}\int_{D}|\nabla \phi|^p dx $$

Next we will prove that 
\begin{eqnarray*}
\frac{1}{p}\int_{w^{\epsilon}>\theta}|\nabla u^{\epsilon}_{-} |^pdx&\geq& -\frac{p-1}{p}\int_{D}|\nabla w^\epsilon \phi|^p dx-\int_{D} |\nabla w^\epsilon \phi|^{p-2}\nabla w^\epsilon \cdot {\nabla u^{\epsilon}_{-}}  \phi dx\\
&-&C\theta-C{\theta}^{\frac{p-1}{p}} \qquad \qquad \qquad \qquad \qquad \qquad \qquad \qquad (**)
\end{eqnarray*}

In fact by Young's inequality we have that
\begin{eqnarray*}
 -\int_{w^{\epsilon}>\theta}|\nabla w^\epsilon \phi|^{p-2}\nabla w^\epsilon \cdot \nabla u^{\epsilon}_{-}  \phi dx&\leq& \frac{p-1}{p}\int_{w^{\epsilon}>\theta}|\nabla w^\epsilon \phi|^p dx\\
&+&\frac{1}{p}\int_{w^{\epsilon}>\theta} |\nabla u^{\epsilon}_{-} |^p 
\end{eqnarray*}

Then by the proof of  Lemma I (a) 
$$\int_{w^{\epsilon}<\theta}|\nabla w^\epsilon \phi|^{p}dx \leq C \theta$$

And by H\"{o}lder inequality
\begin{eqnarray*}
|\int_{w^{\epsilon}<\theta}|\nabla w^\epsilon \phi|^{p-2}\nabla w^\epsilon \cdot {\nabla u^{\epsilon}_{-}}  \phi dx |&\leq& \{\int_{w^{\epsilon}<\theta}|\nabla w^\epsilon \phi|^{p}dx\}^{\frac{p-1}{p}}\\
&\times&\{\int_{D} |\nabla u^{\epsilon}_{-} |^p dx\}^{\frac{1}{p}}\\
&\leq& C\theta^{\frac{p-1}{p}} 
\end{eqnarray*}
Thus
\begin{eqnarray*}
\frac{1}{p}\int_{w^{\epsilon}>\theta}|\nabla u^{\epsilon}_{-} |^pdx&\geq& -\frac{p-1}{p}\int_{D}|\nabla w^\epsilon \phi|^p dx\\
&-&\int_{D} |\nabla w^\epsilon \phi|^{p-2}\nabla w^\epsilon \cdot {\nabla u^{\epsilon}_{-}}  \phi dx-C\theta-C{\theta}^{\frac{p-1}{p}}
\end{eqnarray*}
which concludes $(**)$.

Now we combine the two estimates: $(*)$ and $(**)$ and apply the key Lemma I, then we have the following:
\begin{eqnarray*}
\liminf_{\epsilon\rightarrow0}\frac{1}{p}\int_{D}|\nabla u^{\epsilon}_{-} |^p &\geq&\liminf_{\epsilon\rightarrow0}\{\frac{1}{p}\int_{w^{\epsilon}\leq\theta}|\nabla u^{\epsilon}_{-} |^pdx +\frac{1}{p}\int_{w^{\epsilon}>\theta}|\nabla u^{\epsilon}_{-} |^pdx\}\\
 &\geq&\int_{D} |\nabla \phi|^{p-2}\nabla \phi \cdot \nabla u^{0}_{-}dx - \frac{p-1}{p}\int_{D}|\nabla \phi|^p dx\\
&-& \frac{p-1}{p}\alpha_{0}\int_{D}|\phi|^p dx +  \alpha_{0}\int_{D}|\phi|^{p-2}\phi\cdot u^{0}_{-} dx-C\theta-C{\theta}^{\frac{p-1}{p}}
\end{eqnarray*}
Since $\theta$ is arbitary small, then
\begin{eqnarray*}
\liminf_{\epsilon\rightarrow0}\frac{1}{p}\int_{D}|\nabla u^{\epsilon}_{-} |^p dx &\geq&\int_{D} |\nabla \phi|^{p-2}\nabla \phi \cdot \nabla u^{0}_{-}dx - \frac{p-1}{p}\int_{D}|\nabla \phi|^p dx\\
&-& \frac{p-1}{p}\alpha_{0}\int_{D}|\phi|^p dx +  \alpha_{0}\int_{D}|\phi|^{p-2}\phi\cdot u^{0}_{-} dx
\end{eqnarray*}
Then if we let $\phi=u^{0}_{-}$ (since the test functions are dense in $W^{1,p}(D)$), we have
$$\liminf_{\epsilon\rightarrow0}\int_{D}|\nabla u^{\epsilon}_{-}|^pdx\geq \int_{D}|\nabla u^{0}_{-}|^{p}dx +\alpha_{0}\int_{D}(u^{0}_{-})^p dx.       \qquad \qquad \qquad \Box$$


\vspace*{2cm}

 \begin{thebibliography}{17}
  \bibitem[AB]{} N. Ansini, A. Braides, Asymptotic analysis of periodically-perforated nonlinear media, J. Math. Pures Appl. 81(2002) 439-451.

  \bibitem[C]{} L.A. Caffarelli, Basic Tools for the Regularity Theory of Nonlinear Elliptic Equations I-V , MSRI Talks, 2005.

  \bibitem[CL]{} L.A. Caffarelli, K. Lee, Viscosity method for homogenization of highly oscillating obstacles,  Indiana Univ. Math. J. 57(2008), no. 4, 1715--1741

  \bibitem[CM]{} L.A. Caffarelli, A. Mellet, Random homogenization of an obstacle problem, Ann. I. H. Poincar$\acute{e}$-AN 26(2009), 375-395.

  \bibitem[CSW]{} L.A. Caffarelli, P.E.Souganidis, L. Wang, Homogenization of fully nonlinear, uniformly elliptic and parabolic partial differental equations in stationary ergodic media, Comm.Pure Appl. Math. 58(3)(2005)319-361.

  \bibitem[CV]{} L.A. Caffarelli, A. Vasseur,  The De Giorgi method for nonlocal  fluid dynamics, Preprint, 2009.

  \bibitem[CC]{} L. Carbone, F. Colombini, On convergence of functionals with unilateral constraints, J.Math. Pures Appl.(9)59(4)(1980) 465-500.

  \bibitem[CF1]{} D.Cioranescu, F.Murat,Un terme $\acute{e}$trange venu d'ailleurs, in: Nonlinear Partial Differential Equations and their Applications. Coll$\grave{e}$ge de France Seminar, vol. II, Paris, 1979/1980, in: Res. Notes in Math, vol.60, Pitman, Boston, MA, 1982, PP.98-138, 389-390.

  \bibitem[CF2]{} D.Cioranescu, F.Murat,Un terme $\acute{e}$trange venu d'ailleurs. II, in: Nonlinear Partial Differential Equations and their Applications. Coll$\grave{e}$ge de France Seminar, vol. III, Paris, 1980/1981, in: Res. Notes in Math, vol.70, Pitman, Boston, MA, 1982, PP.154-178, 425-426.

  \bibitem[D]{} G.Dal Maso, Asymptotic behaviour of minimum problems with bilateral obstacles, Ann. Mat. Pura Appl.(4)129(1981)327-366.

  \bibitem[DL]{} G.Dal Maso, P.Longo, $\Gamma$-limits of obstacles, Ann. Mat. Pura Appl.(4)128(1981)1-50.

  \bibitem[DM]{} G.Dal Maso, L.Modica, Nonlinear stochastic homogenization and ergodic theory, J. Reine Angew. Math. 368(1986) 28-42.

   \bibitem[DDL]{} E. De Giorgi, G.Dal Maso, P.Longo, $\Gamma$-limits of obstacles, Atti Accad. Naz. Lincei Rend. Cl. Sci. Fis. Mat. Natur. (8)68(6)(1980)481-487.

    \bibitem[E]{} L. C. Evans, A new proof of local $C^{1,\alpha }$ regularity for solutions of certain degenerate elliptic P.D.E.  J. Differential Equations  45  (1982), no. 3, 356--373.

   \bibitem[GT]{} D. Gilbarg, N.S. Trudinger, Elliptic Partial Differential Equations of Second Order, Second edition, Grundlehren der Mathematischen Wissenschaften 224, Springer-Verlag, Berlin-New York, 1983.
   \bibitem[HL]{} Q. Han, F. Lin, Elliptic partial differential equations. Courant Lecture Notes in Mathematics, 1. New York University, Courant Institute of Mathematical Sciences, New York; American Mathematical Society,Providence, RI,1997


   \bibitem[L]{}G. M. Lieberman, Boundary regularity for solutions of degenerate elliptic equations.  Nonlinear Anal.  12  (1988),  no. 11, 1203--1219.

   \bibitem[MZ]{} J. Mal\'{y}, W. Ziemer, Fine regularity of solutions of elliptic partial differential equations. Mathematical Surveys and Monographs, 51. American Mathematical Society, Providence, RI, 1997.

  \bibitem[W]{} L. Wang, Compactness methods for certain degenerate elliptic equations.  J. Differential Equations  107 (1994),  no. 2, 341--350.

  \bibitem [Z]{} W. Zimmer, Weakly Differentiable Functions,Sobolev spaces and functions of bounded variation. Graduate Texts in Mathematics, 120. Springer-Verlag, New York, 1989.
\end{thebibliography}
\end{document}